\documentclass{article}% Math and Physical Sciences Reference Style

\usepackage{ctex}
\usepackage{amsmath}
\usepackage{amsfonts,amssymb}
\usepackage{extarrows}
\usepackage{geometry}

\usepackage[utf8]{inputenc}
\newtheorem{theorem}{Theorem}[section]
\newtheorem{definition}{Definition}[section]
\newtheorem{lemma}{Lemma}[section]

\newtheorem{fact}{Fact}[section]
\newtheorem{conjecture}{Conjecture}[section]
\def\q{\hfill\rule{1ex}{1ex}}
\usepackage{graphicx}

\usepackage{enumerate}

\usepackage{
	amsmath,			% Math Environments
	amssymb,			% Extended Symbols
	enumerate,		    % Enumerate Environments
	graphicx,			% Include Images
	lastpage,			% R\part{title}eference Lastpage
	multicol,			% Use Multi-columns
	multirow,			% Use Multi-rows
	pifont,			    % For Checkmarks
}

\usepackage[numbers]{natbib}
\begin{document}

% --- PAPER INFO ---

\title{A polynomial time algorithm to find star chromatic index on bounded treewidth graphs with given maximum degree}
\author{{\small\bf Yichen Wang}\thanks{email:  wangyich22@mails.tsinghua.edu.cn
}\quad {\small\bf Mei
Lu}\thanks{email: lumei@tsinghua.edu.cn}\\
{\small Department of Mathematical Sciences, Tsinghua
University, Beijing 100084, China.}}

%\author{Yichen Wang}
\date{2023 3}

\maketitle\baselineskip 16.3pt

\begin{abstract}
% Let $G=(V,E)$ be a graph.
A star edge coloring of a graph $G$ is a proper edge coloring with no 2-colored path or cycle of length four.
The star edge coloring problem is to find an edge coloring of a given graph $G$ with minimum number $k$ of colors such that $G$ admits a star
edge coloring with $k$ colors. This problem is known to be NP-complete.
In this paper, for a bounded treewidth graph with given maximum degree, we show
that it can be solved in polynomial time.
\end{abstract}

%\textbf{AMS classification: }\textit{05C75, 05C65, 05C05}\vskip 0.3cm

\vskip.3cm

% -------------------------------------------------------------
% --------- Here begins INTRODUCTION SECTION ------------------
\section{Introduction}

% Let $G=(V,E)$ be a graph.
A proper edge coloring of a graph $G$ with vertex set $V$ and edge set $E$ is an assignment of colors to the edges of $G$ such that no two adjacent edges have the same color.
Under additional constraints on the proper edge coloring, we have a variety of colorings such as strong edge coloring, vertex distinguishing coloring and so on.
A \textit{star edge coloring} of $G$ is a proper edge coloring  where at least three distinct
colors are used on the edges of every path and cycle of length four, i.e., there is neither
bichromatic path nor cycle of length four.
We say $G$ is \textit{$k$-star-edge-colorable} if $G$ admits a star edge coloring using at most $k$ colors
and the \textit{star chromatic index} $\chi_{st}'(G)$ of $G$  is the minimum number $k$ for which
$G$ has a star edge coloring with $k$ colors.

Star edge coloring is invited by Liu and Deng~\cite{StarEdgeColoringFirstInvited} motivated by the vertex version~\cite{StarVertexColoring,StarVertexColoring2,StarVertexColoring3}. In the same paper, Liu and Deng  presented an upper bound on the star chromatic index of
graphs with maximum degree $\Delta\ge 7$. In \cite{DMS}, Dvo\v r\'{a}k, Mohar and \v S\'{a}mal  presented some upper bounds and lower
bounds on the star chromatic index of complete graphs and subcubic graphs (i.e. with maximum degree at most three).  Some bounds on the star
chromatic index of subcubic outerplanar graphs, trees and outerplanar graphs can be find in \cite{KPR, StarChromaticIndexComplexity, LMS, PV, 3W}.
The relevant research results of star chromatic index can be referred to the survey \cite{SurveyOnStarEdgeColoring}.

In \cite{StarChromaticIndexComplexity}, Lei, Shi and Song showed that it is $NP$-complete to determine whether $\chi_{st}'(G) \le 3$ for an arbitrary graph $G$.
In \cite{ComputeStarChromaticIndexOfTrees}, Omoomi, Roshanbin and Dastjerdi presented a polynomial time algorithm that finds an optimum star edge-coloring
for every tree.

The treewidth of a graph is an important invariant in  graph theory. The concept of treewidth was originally introduced by Bertel\'e and  Brioschi \cite{BB} under the name of dimension. It was later rediscovered by Halin \cite{Hal} in 1976 and by  Robertson and  Seymour \cite{RS84} in 1984.
 The treewidth of a graph gives an indication of how far away the graph is from being
a tree or forest. %The closer the graph is to being a forest, the smaller is its treewidth.
 The treewidth also
is a  parameter that plays a fundamental role in various graph algorithms. It is well-known that many NP-complete problems can be solved in polynomial time on bounded treewidth graphs \cite{EdgeColoringPartialKTree,StrongEdgeColoringPartialKTree,TotalColoringPartialKTree}.
In this paper, for bounded treewidth graphs with given maximum degree, we propose a polynomial time algorithm which can determine the star chromatic  index of $G$. Our main idea comes from \cite{StrongEdgeColoringPartialKTree}.

The rest of this paper is organized as follows. In section 2, we will give some  terminology and notations that we use in this paper. In section 3, we give a polynomial time algorithm to determine the star
chromatic index of bounded treewidth graphs with given maximum degree. Some more discussion will be given in section 4.

% --------------------------------------------------------------
% --------- Here begins PRELIMINARIES AND DEFINITIONS SECTION ------------------
\section{Preliminaries and definitions}

In this section, we give  definitions involving in treewidth.  The treewidth of a graph is defined through
the concept of tree-decomposition. %In the following, we will use $I$ to denote the vertex set of $T$ when $T$ is a tree and we call the vertex of $T$ {\em node}.

\begin{definition}
	A tree-decomposition of a graph $G = (V, E)$ is a pair $(X, T)$,
	where $T(I, F)$ is a tree with vertex set $I$ and edge set $F$, and $X = \{X_i \mid i \in I\}$ is a family of subsets of $V$, one for each node of $T$, such that:
	\begin{itemize}
		\item $\bigcup \limits_{i \in I}X_i = V$.
		\item for each edge $uv \in E$, there exists an $i \in I$ such that $u,v \in X_i$.
		\item for all $i,j,s \in I$, if $j$ is on the path from $i$ to $s$ in $T$, then $X_i \cap X_s \subseteq X_j$.
	\end{itemize}
\end{definition}

Suppose $i$ is a node of $T$, there exists a vertex set $X_i$ corresponds to $i$.
In the following, we also call $X_i$ the node of $T$ for convenience. If $X_i$ is a non-leaf node of $T$, we call it an internal node.
The width of a tree-decomposition $(X, T)$ is $\max \limits_{i\in I} |X_i| - 1$.
The treewidth of a graph $G$ is the minimum treewidth over all possible tree-decompositions of $G$.
The problem of deciding whether a graph has tree-decomposition of treewidth at most $k$ is NP-complete \cite{A} and Bodlaender \cite{SmoothTreeDecomposition2} proved that the problem is fixed-parameter tractable, that is, there is a polynomial $p$ and an algorithm that, for a given graph $G=(V,E)$, computes a tree-decomposition of $G$ of width $k$ in time at most
$2^{p(k)}|V|$.

We say a tree-decomposition $(X, T)$ of treewidth $k$ is \textit{smooth}, if $T$ is a binary tree and:
\begin{itemize}
	\item for all $i \in I$, $|X_i| = k + 1$.
	\item every internal node $X_i$ has two children, and if $X_L,X_R$ are its children, then either $X_L = X_i$ or $X_R = X_i$.
	\item for all $(i,j) \in F$: $k \le |X_i \cap X_j| \le k+1$.
	\item for each edge $uv \in E$, there is at least one leaf $i \in I$, with $u,v \in X_i$.
\end{itemize}

It can be shown that any tree-decomposition of a graph $G$ can be transformed into a smooth tree-decomposition of $G$ with the same treewidth and size $O(n)$ in linear time~\cite{SmoothTreeDecomposition1,SmoothTreeDecomposition2}.

For a given graph $G$ with treewidth at most $k$, let $(X,T)$ be its tree-decomposition.
For each node $X_i$ of $T$, we define $T_i$ to be the subtree of $T$ rooted at $i$.
 Let $X_i$ be a node of $T$.
 If $X_i$ is a leaf, then we define $V_i = X_i$, $E_i = \{uv\in E ~|~ u,v \in X_i\}$; if $X_i$ is an internal node with children $X_L,X_R$, then  we define $V_i = V_L \cup V_R, E_i = E_L \cup E_R$.
The graph $G(V_i, E_i)$ is denoted by $G_i$.

% The idea of our algorithm is that coloring leaf nodes enumeratively and then color internal nodes recursively from its two children.
% When implementing the algorithm,
% We hope that an edge is uniquely colored in one leaf node.
% For this purpose, we need to define another set of edges.
For every edge $e \in E$, there is at least one leaf $i$ of $T$ such that $u,v \in X_i$.
We choose one such $i$ as the \textit{representative} of $e$ and denote it by $rep(e)$.
If $X_i$ is a leaf, we define $E_i'= \{e \in E ~|~ rep(e) \in T_i\}$; if $X_i$ is an internal node with children $X_L, X_R$, we define $E_i' = E_L' \cup E_R'$.
Note that by this definition, for every node $X_i$, $E_i' \subseteq E_i$, and if $X_i$ is an internal node with children $X_L, X_R$, then $E_L' \cap E_R' = \emptyset$.
We denote the graph on vertex set $V_i$ and edge set $E_i'$ by $G_i'$.

In the rest of this paper, when $X_i$ is an internal node in a smooth tree-decomposition, we always use $X_L$ and $X_R$ to represent its two children.
Without loss of generality, we  assume $X_L = X_i$, $X_i - X_R = \{v'\}$ and $X_R - X_i = \{v''\}$. If $e\in E(G)$, we also use $e$ to denote the subset consisting of two vertices of $e$ for short.
Let $\Delta$ be the maximum degree of $G$. The following facts are obvious from definition.%We usually ignore the subscript in the absense of ambiguity.

\begin{fact}\label{fact: v' not in E_R', v'' not in E_L'}
Let $X_i$ be an internal node in a smooth tree-decomposition and  $X_L$ and $X_R$  its two children.
If $ e \in E_R'$ (resp. $ e \in E_L'$), then $v'\notin e$ (resp. $v''\notin e$).
\end{fact}

\begin{fact}\label{fact: intersecting vertex is in X_i}
Let $X_i$ be an internal node in a smooth tree-decomposition and  $X_L$ and $X_R$ its two children.
If $e_1 \in E_L', e_2 \in E_R'$ and $\{v \}= e_1 \cap e_2$, then $v \in X_i-\{v'\}$.
\end{fact}

% ------------------------------------------------------------------------------
% --------- Here begins MAIN ALGORITHM SECTION ------------------

\section{The star edge coloring algorithm}

In this section, we give a polynomial time algorithm to
find the star  chromatic index of bounded treewidth graphs with given maximum degree.  Here we present the main theorem of our paper.

\begin{theorem}\label{thm: main thm: poly algorithm of star edge coloring on partial k-tree and bounded max degree}
	% For every fixed integers $k, \Delta$, there is a deterministic algorithm that, for given graph $G$ on $n$ vertices with treewidth $k$ and maximum degree $\Delta$, and an integer $c$,
	% determines in time $O(nc^{2{(k+1)}^2\Delta^6})$ whether $G$ has a star edge coloring using at most $c$ colors or not, and if so finds such a star edge coloring.

	For everg graph $G$ of order $n$ with treewidth $k$ and maximum degree $\Delta$, and integer $c$,
	there is a deterministic algorithm that determines in time $O(nc^{2{(k+1)}^2\Delta^6})$ whether $G$ has a star edge coloring using at most $c$ colors or not and finds such star edge coloring if it exists.
\end{theorem}

Let $(X,T)$ be a smooth tree-decomposition of $G$ and  $C = \{1,2,\ldots,c\}$  the set of colors.
For a node $X_i$ of $T$, a mapping $f: E_i' \rightarrow C$ is called a \textit{partial coloring} of $G_i'$.
We say a partial coloring is a \textit{proper} edge coloring if no two adjacent edges have the same color.
We say a partial coloring $f$ of $G_i'$ is \textit{valid} if it is a star edge coloring of $G_i'$, that is, it is proper and no path or cycle of length $4$ in $G_i'$ is bicolored.

Consider a partial coloring $f$ of $G_i'$.
For a fixed color pair $(c_1,c_2) \in C^2$, we define the neighbour sets as the following equations.
% Eq. (\ref{eq: neighbour sets definition}).
\begin{equation}
\begin{aligned}
	N^{(1)}_f(c_1,c_2) & = \{ v_0 \in X_i ~|~ \exists \{v_0v_1, v_1v_2\} \subseteq E_i', f(v_0v_1) = c_1, f(v_1v_2) = c_2 \}, \\
	N^{(2)}_f(c_1,c_2) & = \{ v_0 \in X_i ~|~ \exists \{v_0v_1, v_1v_2\} \subseteq E_i', f(v_0v_1) = c_2, f(v_1v_2) = c_1 \}, \\
	N^{(3)}_f(c_1,c_2) & = \{ (v_0, v_2) \in X_i^2 ~|~ \exists \{v_0 v_1 , v_1 v_2\} \subseteq E_i', f(v_0v_1) = c_1, f(v_1v_2) = c_2 \}, \\
	N^{(4)}_f(c_1,c_2) & = \{ v_0 \in X_i ~|~ \exists \{v_0 v_1, v_1 v_2, v_2 v_3\} \subseteq E_i', f(v_0v_1) = c_1, f(v_1v_2) = c_2, f(v_2v_3) = c_1 \}, \\
	N^{(5)}_f(c_1,c_2) & = \{ v_0 \in X_i ~|~ \exists \{v_0 v_1, v_1 v_2, v_2 v_3\} \subseteq E_i', f(v_0v_1) = c_2, f(v_1v_2) = c_1, f(v_2v_3) = c_2 \}, \\
	N^{(6)}_f(c_1,c_2) & = \{ v_0 \in X_i ~|~ \exists v_0v \in E_i', f(v_0v) = c_1 \}, \\
	N^{(7)}_f(c_1,c_2) & = \{ v_0 \in X_i ~|~ \exists v_0v \in E_i', f(v_0v) = c_2 \}, \\
	N^{(8)}_f(c_1,c_2) & = \{ (v_0, v_1) \in X_i^2 ~|~ \exists v_0 v_ 1\in E_i', f(v_0 v_1) = c_1 \}, \\
	N^{(9)}_f(c_1,c_2) & = \{ (v_0, v_1) \in X_i^2 ~|~ \exists v_0 v_1 \in E_i', f(v_0 v_1) = c_2 \}. \\
\end{aligned}
\label{eq: neighbour sets definition}
\end{equation}

Let $S$ be a set. We use $2^{S}$ to represent the power set of $S$.
%All neighbour sets of $(c_1,c_2)$ is a nine-element tuple whose third, eighth and nineth elements are from $2^{X_i^2}$ and the rest are from $2^{X_i}$.
Denote\begin{equation*}
	\mathcal{S}(X_i) = \{ \big( A^{(1)}, \ldots, A^{(9)} \big) ~|~ A^{(j)} \subseteq X_i, j=1,2,4,5,6,7 \mbox{~and~} A^{(j)} \subseteq X_i^2, j=3,8,9\}.
\end{equation*}

In the following, we abbreviate $\big( A^{(1)}, \ldots, A^{(9)} \big)$ to ${\big( A^{(j)} \big)}_{1 \le j \le 9}$.
For a partial coloring $f$ of $G_i'$, define the color class function $CL_f: \mathcal{S}(X_i) \rightarrow 2^{C^2}$ as Eq. (\ref{eq: definition of CL}).
In this case, we say $CL_f$ is the color class function on $X_i$.

\begin{equation}\label{eq: definition of CL}
	\begin{aligned}
		CL_{f}\left( {\big( A^{(j)} \big)}_{1 \le j \le 9}\right ) = \{ (c_1,c_2) \in C^2 ~|~ A^{(j)} =N^{(j)}_f(c_1,c_2), 1 \le j \le 9 \}.
	\end{aligned}
\end{equation}

A color class function $CL_f$ is \textit{active} if and only if there exists a valid partial coloring $g$ such that $CL_{f} = CL_{g}$.

\begin{fact}\label{fact: CL is a partition} Let $f$ be a partial coloring of $G_i'$. Then
		$\{CL_{f}( {\big( A^{(j)} \big)}_{1 \le j \le 9} ) ~|~ {\big( A^{(j)} \big)}_{1 \le j \le 9}  \in \mathcal{S}(X_i)\}$ forms a partition of $C^2$.
\end{fact}

Our idea of developing an algorithm is to calculate all active color class functions on a node $X_i$.
For any partial coloring $f$ of $G_i'$, the color class function is only concerned with the coloring of edges adjacent to a vertex whose distance is at most $3$ to some vertex in $X_i$
(the distance between two vertices is the minimum length of a path connecting them).
Since at most $(k+1)\Delta^3$ vertices have a distance at most $3$ to some vertex in $X_i$, at most ${(k+1)}^{2}\Delta^6$ edges have an effect on the color class function.
Therefore, at most $c^{{(k+1)}^2\Delta^6}$ color class functions use at most $c$ colors.

What we are most concerned about is how to determine whether a color class function is active.
To do this, we first need to analyse the relationship of color class functions on $X_i,X_L$ and $X_R$.

Let $f$ be a partial coloring of $G_i'$. Define $f_L = f~|~G_L'$ and $f_R = f~|~G_R'$.
	For any $(c_1,c_2) \in C^2$ and $1 \le j \le 9$, let $A^{(j)} = N^{(j)}_{f}(c_1,c_2), A_L^{(j)} = N^{(j)}_{f_L}(c_1,c_2)$ and $A_R^{(j)} = N^{(j)}_{f_R}(c_1,c_2)$. Denote
\begin{equation}
	\begin{aligned}
		\bar{A}^{(1)} & = \{ v_0 ~|~ (v_0, v_1) \in A_L^{(8)}, v_1 \in A_R^{(7)}\} \cup \{ v_0 ~|~ (v_0, v_1) \in A_R^{(8)}, v_1 \in A_L^{(7)}  \}, \\
		\bar{A}^{(2)} & = \{ v_0 ~|~ (v_0, v_1) \in A_L^{(9)}, v_1 \in A_R^{(6)}\} \cup \{ v_0 ~|~ (v_0, v_1) \in A_R^{(9)}, v_1 \in A_L^{(6)}  \}, \\
		\bar{A}^{(3)} & = \{ (v_1, v_2) ~|~ \exists (v_1, v_3) \in A_L^{(8)}, (v_3,v_2) \in A_R^{(9)}\} \cup \{ (v_1, v_2) ~|~ \exists  (v_1, v_3) \in A_R^{(8)}, (v_3,v_2) \in A_L^{(9)}\}, \\
		\bar{A}^{(4)} & = \{ v_0 ~|~ \exists (v_0, v_1) \in A_L^{(3)}, v_1 \in A_R^{(6)}\} \cup \{ v_0 ~|~ \exists  (v_0, v_1) \in A_L^{(8)}, (v_1,v_2) \in A_R^{(9)}, v_2 \in A_L^{(6)}\}  \\
		& \cup \{ v_0 ~|~ \exists (v_0, v_1) \in A_L^{(8)}, v_1 \in A_R^{(2)}\} \cup \{ v_0 ~|~ \exists  (v_0, v_1) \in A_R^{(3)}, v_1 \in A_L^{(6)}\}   \\
		& \cup \{ v_0 ~|~ \exists (v_0, v_1) \in A_R^{(8)}, (v_1,v_2) \in A_L^{(9)}, v_2 \in A_R^{(6)}\} \cup \{ v_0 ~|~ \exists  (v_0, v_1) \in A_R^{(8)}, v_1 \in A_L^{(2)}\},   \\
		\bar{A}^{(5)} & = \{ v_0 ~|~ \exists  (v_1, v_0) \in A_L^{(3)}, v_1 \in A_R^{(7)}\} \cup \{ v_0 ~|~ \exists  (v_0, v_1) \in A_L^{(9)}, (v_1,v_2) \in A_R^{(8)}, v_2 \in A_L^{(7)}\}  \\
		& \cup \{ v_0 ~|~ \exists  (v_0, v_1) \in A_L^{(9)}, v_1 \in A_R^{(1)}\} \cup \{ v_0 ~|~ \exists  (v_1, v_0) \in A_R^{(3)}, v_1 \in A_L^{(7)}\}   \\
		& \cup \{ v_0 ~|~ \exists  (v_0, v_1) \in A_R^{(9)}, (v_1,v_2) \in A_L^{(8)}, v_2 \in A_R^{(7)}\} \cup \{ v_0 ~|~ \exists  (v_0, v_1) \in A_R^{(9)}, v_1 \in A_L^{(1)}\},   \\
		\bar{A}^{(j)} & = \emptyset, ~~6 \le j \le 9. \\
	\end{aligned}
	\label{eq: left and right decomposition 2}
	\end{equation}
Then, we have the following lemmas.
\begin{lemma}\label{lemma: neighbour recursion}
For $1 \le j \le 9$,	we have
		\begin{equation}
	\begin{aligned}
		A^{(j)} & = A_L^{(j)} \cup A_R^{(j)} \cup \bar{A}^{(j)} - \{v''\},\\
	\end{aligned}
	\label{eq: left and right decomposition 1}
	\end{equation}	
where $\{v''\}=X_R - X_i$.	
\end{lemma}

\textbf{Proof:} We just show that the result holds for $A^{(1)}$ and the other results can be
showed by the same way.
Fixed color pair $(c_1,c_2)$, we verify $A^{(1)}$ one by one.

  If $v_0 \in A^{(1)}$, then $v_0 \in X_i$ and so $v_0\not=v''$. By the definition, there exists  $\{v_0v_1, v_1v_2\} \subseteq E_i'$ such that $f(v_0v_1) = c_1$ and $f(v_1v_2) = c_2$.
%Considering $v_0v_1, v_1v_2$ are from $G_L'$ or $G_R'$, there are several situations:

\begin{enumerate}
	\item If $v_0v_1 \in E_L'$ and $v_1v_2 \in E_L'$, then $v_0 \in A_L^{(1)}$.
	\item If $v_0v_1 \in E_R'$ and $v_1v_2 \in E_R'$, then $v_0 \in A_R^{(1)}$.
	\item If $v_0v_1 \in E_L'$ and $v_1v_2 \in E_R'$, from Fact~\ref{fact: intersecting vertex is in X_i}, we have $v_1 \in X_i - \{v' \} = X_L - \{v' \}$.
	By the definition, $(v_0,v_1) \in A_L^{(8)}, v_1 \in A_R^{(7)}$ and then we have $v_0 \in \bar{A}^{(1)}$.
	\item If $v_0v_1 \in E_R'$ and $v_1v_2 \in E_L'$, from Fact~\ref{fact: intersecting vertex is in X_i}, we have $v_1 \in X_i - \{v' \} = X_L - \{v' \}$.
	From Fact~\ref{fact: v' not in E_R', v'' not in E_L'}, we have $v_0 \neq v'$.
	By the definition, $(v_0,v_1) \in A_R^{(8)}, v_1 \in A_L^{(7)}$ and then we have $v_0 \in \bar{A}^{(1)}$.
\end{enumerate}

Therefore, we have $v_0 \in A_L^{(1)} \cup A_R^{(1)} \cup \bar{A}^{(1)} - \{v''\}$ which implies $A^{(1)} \subseteq  A_L^{(1)} \cup A_R^{(1)} \cup \bar{A}^{(1)} - \{v''\}$.
Then, we are going to show $ A_L^{(1)} \cup A_R^{(1)} \cup \bar{A}^{(1)} - \{v''\} \subseteq A^{(1)}$.
It is not hard to prove by considering the four situations conversely. The proof is omitted.
\q

\begin{lemma}\label{lemma: relationship of color class}
Assume $f$ is a partial coloring of $G_i'$. Let $CL_{f}, CL_{f_L}$ and $ CL_{f_R}$ be the color class functions on $X_i$, $X_L$ and  $X_R$ respectively.
Then, for all $ {\big( A^{(j)} \big)}_{1 \le j \le 9}  \in \mathcal{S}(X_i)$, we have
\begin{equation}\label{eq: color class CL recursion}
	CL_f\left( {\big( A^{(j)} \big)}_{1 \le j \le 9} \right) = \bigcup \left(CL_{f_L}\left( {\big( A_L^{(j)} \big)}_{1 \le j \le 9} \right) \cap CL_{f_R}\left( {\big( A_R^{(j)} \big)}_{1 \le j \le 9} \right)\right),
\end{equation}
where the union is taken over all $ {\big( A_L^{(j)} \big)}_{1 \le j \le 9} $ and $ {\big( A_R^{(j)} \big)}_{1 \le j \le 9} $ satisfying Eq. (\ref{eq: left and right decomposition 1}).

\end{lemma}

%Suppose $f$ is a partial coloring of $G_i'$, $f_L = f~|~G_L'$ and $f_R = f~|~G_R'$. Let $CL_{f}, CL_{f_L}$ and $CL_{f_R}$ be the color class functions of $f,f_L$ and $f_R$ respectively.
\textbf{Proof:} Given $ {\big( A^{(j)} \big)}_{1 \le j \le 9}  \in \mathcal{S}(X_i)$, assume $(c_1, c_2) \in CL_f\left( {\big( A^{(j)} \big)}_{1 \le j \le 9}\right )$. Then, we have $N_{f}^{(j)}(c_1,c_2)=A^{(j)}$ for all $1\le j\le 9$. Let $ {\big( A_L^{(j)} \big)}_{1 \le j \le 9} = {\big( N_{f_L}^{(j)}(c_1,c_2) \big)}_{1 \le j \le 9}$ and $ {\big( A_R^{(j)} \big)}_{1 \le j \le 9} = {\big( N_{f_R}^{(j)}(c_1,c_2) \big)}_{1 \le j \le 9} $.
By Fact 3.1, $(c_1,c_2)$ belongs to exactly one $CL_{f_L}\left( {\big( A_L^{(j)} \big)}_{1 \le j \le 9} \right)$ and exactly one $CL_{f_R}\left( {\big( A_R^{(j)} \big)}_{1 \le j \le 9} \right)$.
Hence, we have that $(c_1,c_2) \in CL_{f_L}( {\big( A_L^{(j)} \big)}_{1 \le j \le 9} ) \cap CL_{f_R}( {\big( A_R^{(j)} \big)}_{1 \le j \le 9} )$
and by Lemma~\ref{lemma: neighbour recursion}, Eq. (\ref{eq: left and right decomposition 1}) is satisfied.
Therefore, we have $$CL_f\left( {\big( A^{(j)} \big)}_{1 \le j \le 9} \right) \subseteq \bigcup \left(CL_{f_L}\left( {\big( A_L^{(j)} \big)}_{1 \le j \le 9} \right) \cap CL_{f_R}\left( {\big( A_R^{(j)} \big)}_{1 \le j \le 9} \right)\right).$$

Assume $(c_1, c_2) \in CL_{f_L}( {\big( A_L^{(j)} \big)}_{1 \le j \le 9} ) \cap CL_{f_R}( {\big( A_R^{(j)} \big)}_{1 \le j \le 9} )$ for some $  {\big( A_L^{(j)} \big)}_{1 \le j \le 9} ,  {\big( A_R^{(j)} \big)}_{1 \le j \le 9} $ satisfying Eq. (\ref{eq: left and right decomposition 1}). Then, $(c_1, c_2) \in CL_{f_L}( {\big( A_L^{(j)} \big)}_{1 \le j \le 9} )$ and $(c_1, c_2) \in CL_{f_R}( {\big( A_R^{(j)} \big)}_{1 \le j \le 9} )$. So $ {\big( N_{f_L}^{(j)}(c_1,c_2) \big)}_{1 \le j \le 9}={\big( A_L^{(j)} \big)}_{1 \le j \le 9} $ and $ {\big( N_{f_R}^{(j)}(c_1,c_2) \big)}_{1 \le j \le 9} ={\big( A_R^{(j)} \big)}_{1 \le j \le 9} $.
From Lemma~\ref{lemma: neighbour recursion}, we have $(c_1,c_2) \in CL_{f}( {\big( A^{(j)} \big)}_{1 \le j \le 9} )$
which implies $\bigcup (CL_{f_L}( {\big( A_L^{(j)} \big)}_{1 \le j \le 9} ) \cap CL_{f_R}( {\big( A_R^{(j)} \big)}_{1 \le j \le 9} )) \subseteq CL_f( {\big( A^{(j)} \big)}_{1 \le j \le 9} )$.
\q

\begin{figure}[htbp]
     \begin{center}
       \includegraphics[width=0.9\textwidth]{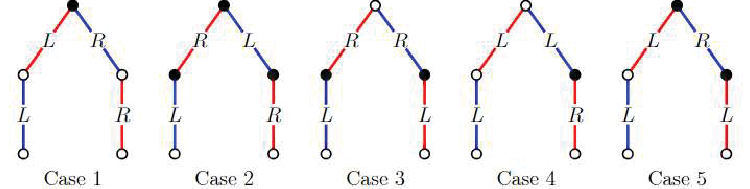}
    \end{center}
        \caption{Five cases of an invalid but proper star coloring in which the bold vertices belong to $X_i - \{v'\}$.
	The two endpoints may be the same and then it results in a bicolored $C_4$ rather than $P_4$.
	Blue and red represent two different colors and $L$ (resp. $R$) represents that the edge is from $G_L'$ (resp. $G_R'$).
	The bold vertices must be in $X_i$ from previous facts.}
    \end{figure}
% FIGURE
%\begin{figure}[t]\label{fig: 5 cases of a invalid star coloring}
   % \centering
  %  \includegraphics[width=\linewidth]{./src/5ContradictionsOfStarColoring.jpg}
  %  \caption{Five cases of an invalid star coloring in which the bold vertices belong to $X_i - \{v'\}$.
	%The two endpoints may be the same and then it results in a bicolored $C_4$ rather than $P_4$. }
%\end{figure}
%\vskip.3cm
Lemma~\ref{lemma: relationship of color class} describes the relationship of color class functions of $X_i, X_L$ and $X_R$.
The remaining problem is how to justify whether a color class function is active.
Lemma~\ref{lemma: active color class justification} provides a method.

\begin{lemma}\label{lemma: active color class justification}
	A color class function $CL_f$  on $X_i$ is active
	if and only if there exists active color class functions $CL_{g_L}$ and $CL_{g_R}$ on $X_L$ and $X_R$ respectively such that for all $ {\big( A^{(j)} \big)}_{1 \le j \le 9}  \in \mathcal{S}(X_i)$,%the following condition holds:
	\begin{equation}\label{eq: color class CL recursion 2}
		CL_f\left( {\big( A^{(j)} \big)}_{1 \le j \le 9} \right) = \bigcup \left(CL_{g_L}( {\big( A_L^{(j)} \big)}_{1 \le j \le 9} ) \cap CL_{g_R}( {\big( A_R^{(j)} \big)}_{1 \le j \le 9} )\right),
	\end{equation}
where the union is taken over all $ {\big( A_L^{(j)} \big)}_{1 \le j \le 9} $ and $ {\big( A_R^{(j)} \big)}_{1 \le j \le 9} $ satisfying Eq. (\ref{eq: left and right decomposition 1}).

	Moreover, if $CL_{g_L}\left( {\big( A_L^{(j)} \big)}_{1 \le j \le 9} \right) \cap CL_{g_R}\left( {\big( A_R^{(j)} \big)}_{1 \le j \le 9}\right) \neq \emptyset$ for some $ {\big( A_L^{(j)} \big)}_{1 \le j \le 9} \in \mathcal{S}(X_L)$ and $ {\big( A_R^{(j)} \big)}_{1 \le j \le 9} \in \mathcal{S}(X_R)$, then the following results hold.
	\begin{enumerate}[(i)]
		\item $A_L^{(1)} \cap A_R^{(2)} = A_L^{(2)} \cap A_R^{(1)} = \emptyset$   (prevent invalid case 1 in Fig. 1).
		\item If $(v_0, v_1) \in A_L^{(8)}, (v_1, v_2) \in A_R^{(9)}$, then $v_0 \notin A_R^{(7)}$ or $v_2 \notin A_L^{(6)}$;
				if $(v_0, v_1) \in A_R^{(8)}, (v_1, v_2) \in A_L^{(9)}$, then $v_0 \notin A_L^{(7)}$ or $v_2 \notin A_R^{(6)}$    (prevent invalid case 2 in Fig. 1).
		\item If $(v_1, v_2) \in A_L^{(3)}$, then $v_1 \notin  A_R^{(7)}$ or $v_2 \notin A_R^{(6)}$;
				if $(v_1, v_2) \in A_R^{(3)}$, then $v_1 \notin  A_L^{(7)}$ or $v_2 \notin A_L^{(6)}$ (prevent invalid case 3 in Fig. 1).
		\item $ A_L^{(4)} \cap  A_R^{(7)} =  A_R^{(4)} \cap  A_L^{(7)} = \emptyset$ (prevent invalid case 4 in Fig. 1).
		\item If $(v_1, v_2) \in  A_L^{(8)}$, then $v_1 \notin  A_R^{(2)}$ or $v_2 \notin A_R^{(7)}$;
				if $(v_1, v_2) \in  A_R^{(8)}$, then $v_1 \notin  A_L^{(2)}$ or $v_2 \notin A_L^{(7)}$ (prevent invalid case 5 in Fig. 1).
		\item  $A_L^{(6)} \cap A_R^{(6)} = A_L^{(7)} \cap A_R^{(7)} = \emptyset$ (assure a proper edge coloring).
	\end{enumerate}
\end{lemma}

\textbf{Proof}
\textit{Necessity:}
Assume that the color class function $CL_f$ of $X_i$ is active and $f$ is a valid partial coloring of $G_i'$.
It is obvious $CL_{f_L}$ defined by $f_L$ on $X_L$ and $CL_{f_R}$ defined by $f_R$ on $X_R$ are active color class functions.
Let $g_L = f_L$ and $g_R = f_R$. By Lemma \ref{lemma: relationship of color class}, Eq.~(\ref{eq: color class CL recursion 2}) holds. Suppose there is $(c_1,c_2)\in CL_{g_L}( {\big( A_L^{(j)} \big)}_{1 \le j \le 9} ) \cap CL_{g_R}( {\big( A_R^{(j)} \big)}_{1 \le j \le 9})$ for some $ {\big( A_L^{(j)} \big)}_{1 \le j \le 9} \in \mathcal{S}(X_L)$ and $ {\big( A_R^{(j)} \big)}_{1 \le j \le 9} \in \mathcal{S}(X_R)$. Then, $(c_1,c_2)\in CL_{f_L}( {\big( A_L^{(j)} \big)}_{1 \le j \le 9} )$ and $(c_1,c_2) \in CL_{f_R}( {\big( A_R^{(j)} \big)}_{1 \le j \le 9})$.
We are going to show that ($i$) to ($vi$) hold. Since $f$ is a valid partial coloring of $G_i'$, ($vi$) holds obviously.

($i$) Assume without loss of generality $v_0\in A_L^{(1)} \cap A_R^{(2)}$. Then there are $v_1,v_2\in X_L$ (resp. $v_3,v_4\in X_R$ such that $f(v_1 v_0) = c_1, f(v_1v_2) = c_2$ (resp. $f(v_3 v_0) = c_2, f(v_3v_4) = c_1$).
Thus we have a bicolored $P_4$ or $C_4$, a contradiction with $f$ being a valid partial coloring of $G_i'$. Hence, ($i$) holds.

The proofs of ($ii$) to ($v$) are the same. We omit here.
%If condition (1) is violated, suppose there exists a color pair $v \in A^{(1)} \cap B^{(2)}$.
%Then, it forms a bicolored $P_4$ or $C_4$ which contradicts the valid hypothesis.
%Similarly we can get such contradictions if $ A_L^{(2)} \cap A_R^{(1)} \neq \emptyset$. Condition (1) actually prevent the first invalid structures in Fig. (\ref{fig: 5 cases of a invalid star coloring}) from existing.

%Similarly condition (2,3,4,5) exactly prevent the last four invalid structures in Fig. (\ref{fig: 5 cases of a invalid star coloring}) from existing and condition (6) asures it is a proper edge coloring.
%The complete proof is omitted here.

\textit{Sufficiency:}
Suppose there exist active color class functions $CL_{g_L}$ and $CL_{g_R}$  on $E_L'$ and $ E_R'$ respectively.
Define $f' : E_i' \rightarrow C$ as following: $f'(e) = g_L(e)$ if $e \in E_L'$ and $f'(e) = g_R(e)$ if $e \in E_R'$.

From Lemma~\ref{lemma: relationship of color class}, we have $CL_{f'} = CL_{f}$. In order to show that $CL_f$  on $X_i$ is active, we just need to prove $f'$ is a valid partial coloring by the definition. By ($vi$), and $CL_{g_L}$ and $CL_{g_R}$ are active on $X_L$ and $X_R$ respectively, we have $f'$ is  a proper edge coloring.
Hence $f'_L$ and $f'_R$ are valid colorings on $X_L$ and $X_R$.
% ------- TBD --------- 他说有问题

Suppose there is a bicolored $P_4$ or $C_4$ on $G_i$ after the coloring $f'$. Then, the $P_4$ or $C_4$ can not belong to only $G_L'$ or $G_R'$.
So it must corresponds to one of the cases in Fig. 1.

Suppose the bicolored $P_4$ or $C_4$ is case 1 in Fig. 1.  Let $E(P_4)$ or $E(C_4)$ = $\{v_1v_2, v_2v_3, v_3v_4, v_4v_5\}$ (when $v_1=v_5$, it is the $C_4$) where $v_1v_2, v_2v_3 \in E_L', v_3v_4, v_4v_5 \in E_R'$ and $f'(v_1v_2) = f'(v_3v_4) = c_1, f'(v_2v_3) = f'(v_4v_5) = c_2$.
From Fact~\ref{fact: intersecting vertex is in X_i}, we have $v_3 \in X_i$.
Then there exist $ {\big( A_L^{(j)} \big)}_{1 \le j \le 9}  \in \mathcal{S}(X_L)$ and $ {\big( A_R^{(j)} \big)}_{1 \le j \le 9}  \in \mathcal{S}(X_R)$ satisfying $(c_1,c_2) \in CL_{g_L}( {\big( A_L^{(j)} \big)}_{1 \le j \le 9} ) \cap CL_{g_R}( {\big( A_R^{(j)} \big)}_{1 \le j \le 9} )$.
By the definition of neighbour sets, we have $v_3 \in A_L^{(1)} \cap B_R^{(2)}$, a contradiction with ($i$).
Similar result can be gained if $v_1v_2, v_2v_3 \in E_R', v_3v_4, v_4v_5 \in E_L'$.
%It shows if case 1 in Fig.~\ref{fig: 5 cases of a invalid star coloring} exists, condition (1) is violated.

By the same argument,  cases 2 to 4 in Fig. 1 can not occur by the conditions ($ii$) to ($v$).
Thus $f'$ is a valid partial coloring which implies $CL_{f}$ is active. \q

\vspace{2em}
\noindent{\bf Proof of Theorem \ref{thm: main thm: poly algorithm of star edge coloring on partial k-tree and bounded max degree}}
In order to prove Theorem \ref{thm: main thm: poly algorithm of star edge coloring on partial k-tree and bounded max degree}, we describe our algorithm.
For a given graph $G$ and a number of colors $c$, we first calculate the smooth tree-decomposition of $G$.
Then, we calculate every active color class functions on each node $X_i$ in the following way.

\begin{enumerate}
	\item  $X_i$ is a leaf node. Enumerate all active $CL_f$ of some $f$ on $X_i$.
	There are at most $c^{{(k+1)}^2\Delta^6}$ color class functions.
	\item  $X_i$ is an internal node.
	For each active color class functions $CL_{f_L}$ and $CL_{f_R}$ on $X_L$ and $X_R$ respectively,
	we calculate all active color class function $CL_f$ of $f$ on $X_i$ by Lemma~\ref{lemma: active color class justification}.
	There are at most $c^{2{(k+1)}^2\Delta^6}$ pair of active $CL_{f_L}, CL_{f_R}$ and verify Lemma~\ref{lemma: active color class justification} for each color class pair takes $O(1)$ time.
\end{enumerate}

Suppose $X_{root}$ is the root node in the smooth tree-decomposition.
Finally, there is an active color class function $CL_{f_r}$ of $X_{root}$ if and only if $G$ have a valid star edge coloring within $c$ colors.
Since there are at most $O(n)$ nodes in the smooth tree-decomposition, the whole algorithm takes $O(nc^{2{(k+1)}^2\Delta^6})$ time.

From the above Lemmas, we can easily prove that there is a valid partial of $X_i$ if and only if the corresponding color class function $CL$ can be computed by the above algorithm.
We did not describe how to find a star edge coloring using the specified number of colors.
But this can be done easily, by doing slight modifications in the algorithm (using extra bookkeeping to store coloring information). We skip the details.

This completes the proof of Theorem~\ref{thm: main thm: poly algorithm of star edge coloring on partial k-tree and bounded max degree}. \q

\vskip.2cm
It is worth to mention that the algorithm can run parallelly in $O(c^{2{(k+1)}^2 \Delta^6}\log n )$ time on $O(n)$ terminals.

% -----------------------------------------------------------
% --------- Here begins DISCUSSION SECTION ------------------

\section{Discussion}

In our algorithm, $\Delta$ must be considered a constant because we need to count the number of color class functions.
To remove the restriction of $\Delta$, we define color class count function $CF\left( {\big( A^{(j)} \big)}_{1 \le j \le 9}\right ) = \left |CL\left( {\big( A^{(j)} \big)}_{1 \le j \le 9} \right)\right|$.
Since $CL$ is a partition of $C^2$, we have:

\begin{equation*}
	\sum_{ {\big( A^{(j)} \big)}_{1 \le j \le 9} \in \mathcal{S}(X_i) } CF\left( {\big( A^{(j)} \big)}_{1 \le j \le 9} \right) = c^2.
\end{equation*}

The number of color class count functions is at most $c^{2^{3k^2+6k}}$.
If the following Conjecture~\ref{conj: color count function equivalence} is proved to be true, we can store $CF$ for each $X_i$ rather than $CL$, and design a polynomial time algorithm.

\begin{conjecture}\label{conj: color count function equivalence}
If $f,g$ are two partial coloring of $G_i'$ with the same color class count function, i.e. $CF_{f}\left( {\big( A^{(j)} \big)}_{1 \le j \le 9} \right) = CF_{g}\left( {\big( A^{(j)} \big)}_{1 \le j \le 9}\right )$ for all $ {\big( A^{(j)} \big)}_{1 \le j \le 9}  \in \mathcal{S}(X_i)$,
then there exists a permutation $\pi$ of colors such that $CL_{\pi \circ f}\left( {\big( A^{(j)} \big)}_{1 \le j \le 9}\right ) = CL_{g}\left( {\big( A^{(j)} \big)}_{1 \le j \le 9} \right)$, for all $ {\big( A^{(j)} \big)}_{1 \le j \le 9}  \in \mathcal{S}(X_i)$.
\end{conjecture}

% ----------------------
% References

\section*{Acknowledgement}
This research was supported by the National Natural Science Foundation of China (Grant 12171272 \&12161141003).

% \bibliography{ref.bib}
% \bibliographystyle{wyc}
% \printbibliography

\end{document}